\documentclass{article}
\usepackage{amsmath,amssymb,amsthm}

\newtheorem{theorem}{Theorem}[section]
\newtheorem{corollary}[theorem]{Corollary}

\newtheorem{lemma}[theorem]{Lemma}
\theoremstyle{remark}
\newtheorem{remark}{Remark}[section]
\theoremstyle{definition}
\newtheorem{definition}{Definition}[section]

\theoremstyle{definition}

\newcommand{\bkR}{\mathbb{R}}
\newcommand{\Ci}{\ensuremath{C^1([a,b];E)}}
\newcommand{\Cii}{\ensuremath{C^2([a,b];E)}}
\newcommand{\Luu}[2]{\mathcal{L}(#1,#2)}

\newcommand{\LFE}{\Luu{F}{E}}
\newcommand{\Auu}[2]{\overrightarrow{(#1,#2)}}

\newcommand{\AFE}{\Auu{F}{E}}
\newcommand{\NLuu}[2]{N_{\Luu{#1}{#2}}}
\newcommand{\NL}{\NLuu{F}{E}}

\newcommand{\Fu}[1]{\overrightarrow{{#1}^\ast}}
\newcommand{\Fast}{\Fu{F}}

\begin{document}

\title{First Integrals for Problems of Calculus\\
of Variations on Locally Convex Spaces\thanks{Presented at
OTFUSA'2005, International Conference on \emph{Operator Theory,
Function Spaces and Applications}, dedicated to the 60th birthday
of Professor F.-O. Speck, 7--9 July 2005, Aveiro, Portugal.
Research Report CM05/I-32.}}

\author{Eug\'{e}nio~A.~M.~Rocha \\ \texttt{eugenio@mat.ua.pt} \\
        \and
        Delfim~F.~M.~Torres \\ \texttt{delfim@mat.ua.pt}}
\date{ Control Theory Group (\textsf{cotg})\\
       Department of Mathematics\\
       University of Aveiro\\
       3810-193 Aveiro, Portugal}

\maketitle

\begin{abstract}
The fundamental problem of calculus of variations is considered
when solutions are differentiable curves on locally convex spaces.
Such problems admit an extension of the Euler-Lagrange equations
[Orlov 2002] for continuously normally differentiable Lagrangians.
Here, we formulate a Legendre condition and an extension of the
classical theorem of Emmy Noether, thus obtaining first integrals
for problems of the calculus of variations on locally convex
spaces.
\end{abstract}

\smallskip


\smallskip


\smallskip



\numberwithin{equation}{section}


\section{Introduction}

The fundamental problem of calculus of variations (CV) is studied
in the setting of infinite dimensional differential geometry
\cite{KMbook}, i.e. where solutions are differentiable curves on
locally convex spaces. The usual problem of CV concern to find
among all functions with prescribed boundary condition, those
which minimize a given functional, i.e.
\begin{eqnarray*}
\int_\Omega f(x,u(x),\nabla u(x))\,dx\longrightarrow\min\, \\
\mbox{s.t.}\:\:\:\:\:\:\:u\in X\;\;\mbox{ and }\:\: u|_{\delta\Omega}=u_0,
\end{eqnarray*}
where $\Omega\subset\bkR^n$ is a bounded open set,
$u:\Omega\subset\bkR^n\rightarrow\bkR^m$, $\nabla u\in\bkR^{nm}$,
$f:\Omega\times\bkR^m\times\bkR^{nm}\rightarrow\bkR$ is a
continuous function, $u_0$ is a given function, and $X$ is a
Banach space. It is well known that problems of CV have very wide
applications in several fields of mathematics, and in many areas
of physics, economics, and biology. In recent years, part of the
renewal of interest in variational methods finds its origins in
nonlinear elasticity \cite{Dacorogna}.

The present work deal with an extension of the setting of the
previous problem by replacing $\bkR^m$ by a locally convex space
$E$. Although, for technical reasons, we will only deal with the
case $n=1$. The problem of the CV on a locally convex space $E$ is
then
\begin{equation}
\label{eq:probCVar}
\textsl{J}\left[x\right] =
\int_a^b L\left(t,x(t),\dot{x}(t) \right)dt \longrightarrow
\min\, ,
\end{equation}
where $L:[a,b]\times E\times E\rightarrow \bkR$,
$x:[a,b]\rightarrow E$. But it is not completely defined without
introducing the precise notions of differentiability and mapping
regularity. It is well known that such functionals
(\ref{eq:probCVar}) are not, in general, Fr\'echet differentiable
(see \textrm{e.g.} \cite{Curtain,Or11}). Comparing with the
classical CV, the main difficulties arise from the substitution of
$\bkR^m$ by a locally convex space $E$; moving from a
finite-dimension vector space with an inner product to an
infinite-dimension vector space with only a family of semi-norms.
The motivation for such problem is the first author interest
\cite{Rocha} on studying calculus of variation and control theory
on the setting of infinite dimensional differential geometry
(differential calculus of smooth mappings between subsets of
sequentially complete locally convex spaces) developed by
Fr\"olicher, Kriegl and Michor \cite{KMbook}; so called convenient
spaces.

A central result of CV is the first order necessary optimality
condition, asserting that optimal solutions satisfy the
Euler-Lagrange equation. Solutions of the Euler-Lagrange equations
are called \emph{extremals}. Extremals include optimal solutions
but, in general, may also include non-optimal solutions. A
function $C(t,x,v)$ which is preserved along all the extremals
(\textrm{i.e.} $C(t,x(t),\dot{x}(t))$ is constant $\forall
t\in[a,b]$ and for any extremal $x$) is called a \emph{first
integral}. The equation $C(t,x(t),\dot{x}(t)) = \text{constant}$
is the corresponding \emph{conservation law}. Conservation laws
are a useful tool to simplify the problem of finding minimizers
\cite{Brunt,Logan}. Emmy Noether was the first to establish a
general theory of conservation laws in the calculus of variations
\cite{MR0406752}. Noether's theorem comprises a universal
principle, connecting the existence of a group of transformations
under which the functional to be minimized is invariant (the
existence of variational symmetries) with the existence of
conservation laws. Noetherian conservation laws play an important
role on a vast number of disciplines, ranging from classical
mechanics, where they find important interpretations such as
conservation of energy, conservation of momentum, or conservation
of angular momentum, to engineering, economics, control theory and
their applications \cite{Gugushvili}.

In this work we use the notion of compactly normally differential
functionals (introduced in \cite{Orlov02}) and inductive scales
of locally convex spaces: (1) to formulate a Legendre condition;
(2) to extend the classical Noether's theorem to the calculus of
variations on locally convex spaces.

The use of inductive scales of locally convex spaces is not a
merely generalization, it is a need. First, as shown in
\cite{Omori78}, Banach manifolds are not suitable for many
questions of global analysis. Second, we require the evaluation
map $E^\ast\times E\rightarrow\bkR$ to be jointly continuous in
order to be able to use integration by parts on some working
space. However, if $E$ is a locally convex space, and $E^\ast$ is
its dual equipped with any locally convex topology, then the
jointly continuity of the evaluation map imply that, in fact, $E$
is a normable space. Since, then, there are neighborhoods
$U\subset E$ and $V\subset E^\ast$ of zero such that the image of
$V\times U$ by the evaluation map is contained on $[-1,1]$. But
then $U$ is contained in the polar of $V$, so it is bounded in
$E$. Therefore, $E$ admits a bounded neighborhood.


\section{Inductive scales of locally convex spaces}

In what follows, $F,E$ are locally convex spaces (LCSs) with the
corresponding determining systems of semi-norms
$\{\|\cdot\|_p\}_{p\in S_F}, \{\|\cdot\|^q\}_{q\in S_E}$, that are
inductively ordered according to the increase of semi-norms. The
set of linear continuous maps from $F$ to $E$ will be denoted by
$\LFE$. For a gentle introduction to locally convex topological
vector spaces we refer the reader to \cite{Curtain}.

Let $A\in\LFE$. For any $q\in S_E$, the normal index of $A$ is the
increasing multi-valued mapping $n_A:S_E\rightarrow 2^{S_F}$
defined by
\begin{equation*}
 n_A(q)=\left\{p\in S_F : \sup_{\|y\|_p\leq 1} \|Ay\|^q
 <+\infty\right\},
\end{equation*}
and $\NL=\{ n_A : A\in\LFE \}$ is the set of all normal indices.
We will consider the following inductive scale of LCSs \cite{Or4}
\begin{equation*}
  \AFE = \{(X_n,\tau_n)\}_{n\in\NL}\:\:\:\mbox{ with
  }\:\:X_n=\{A\in\LFE : n_A\leq n\},
\end{equation*}
\textrm{i.e.} a system of LCSs inductively ordered according to the
continuous embedding $m\leq n\Rightarrow X_{m}\subseteq
X_{n}$; where each space $X_n$ has the projective topology
$\tau_n$ with respect to the determining system of semi-norms
\begin{equation*}
\|A\|_p^q = \left\{\sup_{\|y\|_p\leq 1}\|A y\|^q : p\in n(q), q\in
S_E\right\}.
\end{equation*}
This inductive scale of LCSs generalize classical interpolation
scales \cite{Or6}. Properties of a scale are related with
properties of the spaces of the scale and vice-versa \cite{Or1}.
Convergence in the scale $\AFE$ is the convergence in any space
$X_n$ of the scale. For $Z$ a LCS, an operator $A\in\Luu{Z}{\AFE}$
if $A\in\Luu{Z}{X_n}$ for some $X_n$ and $n\in\NL$.
\begin{definition}
A mapping $\phi:F\rightarrow\AFE$ is
\begin{enumerate}
 \item \emph{continuous} at a point $y_0\in F$ if $y\rightarrow y_0\Rightarrow
 \phi(y)\rightarrow\phi(y_0)$ for some $X_n$.
 \item \emph{uniformly continuous} on a set $D\subset F$ if, for some $X_n$, the map
 $\bar{\phi}:D\rightarrow X_n$ is
 uniformly continuous.
\end{enumerate}
\end{definition}
Let $F_1,F_2$ be LCSs. The canonical isomorphism \cite{Or3}
correspondence $B(y_1,y_2)=(Ay_1)y_2$ between linear operators
$A:F_1\rightarrow\Luu{F_2}{E}$ and bilinear operators
$B:F_1\times F_2\rightarrow E$ justifies the following isometrically identification
\begin{equation}
\label{identif}
  \Auu{F_1}{\Auu{F_2}{E}} \cong \Auu{F_1\times F_2}{E}.
\end{equation}
If $E\equiv\bkR$, the conjugate space has a normal decomposition
into the following inductive scale of Banach spaces
\begin{equation*}
 \Auu{F}{\bkR}\equiv\Fast =\{f\in \Luu{F}{\bkR}:
 \|f\|^p \equiv\sup_{\|x\|_p\leq 1} |f(x)|<+\infty\mbox{ and }p\in S_F\}.
\end{equation*}
In this case, we have the following canonical isometrical isomorphism
of linear and bilinear operators in LCSs
\begin{equation*}
  \Auu{F_1}{\Fu{F_2}} \cong \Fu{(F_1\times F_2)}.
\end{equation*}


\section{Compactly normal differentiability}

Let $F,E$ be LCSs, $y\in F$, and $C$ a convex compact subset of
$F$ having $y$ as limiting point.
\begin{definition}
A map $g:F\rightarrow E$ is
\begin{enumerate}
\item \emph{normally differentiable} at the
point $y$ if $\Delta g(y,h)=g'(y) h+\phi_y(h)$ where
$g'(y)\in\mathcal{L}(F,E)$ and
\begin{equation*}
 \forall q\in S_E\:\exists p\in S_F\: :\:
\frac{\|\phi_y(h)\|^q}{\|h\|_p}\:\stackrel{\:\:h\rightarrow 0\:\:}
{\overrightarrow{\:\:\:\:\:\:\:\:\:\:\:\:\:\:}}\: 0.
\end{equation*}
\item \emph{continuously normally differentiable} at $y$ if $g$ is
normally differentiable in a neighborhood of $y$ and the derived
mapping $g':F\rightarrow\AFE$ is continuous at $y$. This last
condition means that $\tilde{g}':[\alpha,\beta]\subset
F\rightarrow X_n$ is continuous for some $X_n$, and that the
compactness of $[\alpha,\beta]$ implies $n_{\tilde{g}'(y)}\leq n$
for $n\in\NL$ and all $y\in [\alpha,\beta]$. \item \emph{twice
continuously normally differentiable} at $y$ if $g$ is
continuously normally differentiable in a neighborhood of $y$ and
$g':F\rightarrow X_n$ is normally differentiable at a vicinity of
$y$ for some $X_n$. In this case, $g''(y)$ will denote $(g')'(y)$
and, by the identification (\ref{identif}),
$f'':F\rightarrow\Auu{F\times F}{E}$. \item
\emph{K-differentiable} (compactly normally differentiable) at a
point $y$ if for each $C$ the restriction $f=g|_C$ is normally
differentiable at the point $y$. The value $(f|_C)'(y)$ does not
depend on the choice of subset $C$, and it is denoted by
$g'_K(y)$.
\end{enumerate}
\end{definition}

\begin{lemma}\label{fund:lemma} A mapping $g:F\rightarrow E$ is continuously
normally differentiable at a convex compact set $C\subset F$ if
and only if
\[ \frac{\|g(x+h)-g(x)-g'(x)h\|^s}{\|h\|_m}\rightarrow 0\mbox{ as }
h\rightarrow 0, x\in C, \] for some normal index $n_A\in\NL$ and
any $s\in S_E$ and $m\in n_A(s)$.
\end{lemma}

\begin{remark}
If $F$ is a Banach space then $\phi_y(h)=o(\|h\|)$ and $g'$ is the
Fr\'echet derivative of $g$. In such case, we denote the
derivative by a dot, $\dot{g}$.
\end{remark}

We denote by $\Ci$ the space of continuous differentiable mappings
$x:[a,b]\rightarrow E$ with a determining system of
semi-norms $\{\|x\|_1^p\}_{p\in P}$ where
\begin{equation*}
\|x\|_1^p=\sup_{a\leq t\leq b}\|x(t)\|_p+\sup_{a\leq t\leq
b}\|\dot{x}(t)\|_p.
\end{equation*}


\section{Euler-Lagrange and Legendre conditions}

Let $F\equiv [a,b]\times E\times E$. The following theorem proves
that optimal solutions of the problem of the Calculus of
Variations, for $K$-differential functionals, verify an
Euler-Lagrange equation.

\begin{theorem}[\emph{\cite{Orlov02}}]
\label{teo-orlov}
 Let a function $L(t,x,v)$ be continuously normally differentiable
 on $[a,b]\times E\times E$. If the functional $\textsl{J}[\cdot]$
 (\ref{eq:probCVar}) has an extremum at a point
 $x\in\Ci$, then $\textsl{J}[\cdot]$ is $K$-differentiable at $x$ and we have
 \begin{equation}
 \label{eq:teo-orlov}
   \textsl{J}^{\,\prime}_K[x]h=\int_a^b\left[
    \frac{\partial L}{\partial x}(t,x(t),\dot{x}(t))h(t)
     +\frac{\partial L}{\partial v}(t,x(t),\dot{x}(t))\dot{h}(t)
   \right]\,dt =0.
 \end{equation}
\end{theorem}

Observe that $\dot{x}$ is the Fr\'echet derivative of $x$,
where $[a,b]\rightarrow\Auu{[a,b]}{E}$ is identified with
$[a,b]\rightarrow E$. For a given $x$, let $\mu_x:[a,b]\rightarrow F$
be defined by $\mu_x(t)=(t,x(t),\dot{x}(t))$.
For any $(\bar{t},\bar{x},\bar{v})\in F$
define $L^t_{(\bar{x},\bar{v})}:\bkR\rightarrow\bkR$
as $L^t_{(\bar{x},\bar{v})}(t)=L(t,\bar{x},\bar{v})$,
$L^x_{(\bar{t},\bar{v})}:E\rightarrow\bkR$
as $L^x_{(\bar{t},\bar{v})}(x)=L(\bar{t},x,\bar{v})$,
and $L^v_{(\bar{t},\bar{x})}:E\rightarrow\bkR$
as $L^v_{(\bar{t},\bar{x})}(v)=L(\bar{t},\bar{x},v)$.
Partial derivatives are defined as usual
\begin{equation}
\label{eq:partDer}
 \frac{\partial L}{\partial t}(\bar{t},\bar{x},\bar{v})
   =(L^t_{(\bar{x},\bar{v})})'(\bar{t}) \, , \,
 \frac{\partial L}{\partial x}(\bar{t},\bar{x},\bar{v})
   =(L^x_{(\bar{t},\bar{v})})'(\bar{x}) \, , \,
 \frac{\partial L}{\partial v}(\bar{t},\bar{x},\bar{v})
   =(L^v_{(\bar{t},\bar{x})})'(\bar{v}) \, ,
\end{equation}
hence, for a given extremal $\hat{x}$, we have
$\frac{\partial L}{\partial x}\circ\mu_{\hat{x}}:[a,b]\rightarrow\Fu{E}$ and
$\frac{\partial L}{\partial v}\circ\mu_{\hat{x}}:[a,b]\rightarrow\Fu{E}$.
Now, by virtue of the jointly continuous of the evaluation map on $\Fu{E}\times{E}$,
the space $\mathcal{C}\equiv\{A\otimes B\in\Luu{[a,b]}{\Fu{E}}\times\Luu{[a,b]}{E}:
A\mbox{ and }B\mbox{ are differentiable mappings}\}$ is a derivation algebra,
\textrm{i.e.} it admits the Leibniz product rule and the usual integration by parts.

The Euler-Lagrange equation is obtained as a corollary of
Theorem~\ref{teo-orlov} in the usual way, considering the
integration by parts of the second addend of (\ref{eq:teo-orlov})
\begin{equation*}
\int_a^b \left(\frac{\partial L}{\partial v}\circ\mu_{\hat{x}}(t)\right)\dot{h}(t)\,dt
=\left(\left.\frac{\partial L}{\partial v}\circ\mu_{\hat{x}}(t)\right)h(t)\right|_a^b
-\int_a^b \frac{d}{dt}
\left[ \frac{\partial L}{\partial v}\circ\mu_{\hat{x}}(t)\right] h(t)\, dt \, ,
\end{equation*}
and using an extension of the fundamental lemma
of the calculus of variations.

\begin{corollary}[\cite{Orlov02}]
 Let a function $L(t,x,v)$ be twice continuously normally differentiable
 on $[a,b]\times E\times E$. If functional $\textsl{J}[\cdot]$
 (\ref{eq:probCVar}) has an extremum at a point
 $\hat{x}\in\Ci$, then $\hat{x}$ satisfies the Euler-Lagrange equation (in $\Fu{E}$)
 \begin{equation}
   \label{eq:EL:LCS}
   \frac{\partial L}{\partial x}\circ\mu_{\hat{x}}(t)
    -\frac{d}{dt}\frac{\partial L}{\partial v}\circ\mu_{\hat{x}}(t)=0.
 \end{equation}
\end{corollary}
Following the classical terminology \cite{Brunt,Logan},
we call \emph{extremals} to the solutions of \eqref{eq:EL:LCS};
\emph{first integrals} to functions which are kept constant
along all the extremals.

Other necessary conditions exist apart from the Euler-Lagrange
equation. In what follows we obtain, from the second variation,
the so called Legendre condition. Let us observe that the
Euler-Lagrange equation is just a consequence of the Fermat lemma
\emph{(if a function $f:E\rightarrow\bkR$ has a local extremum at
a point $x$ and is normally differentiable at this point, then we
have $f'(x)=0$)}, and the Legendre condition is a consequence of
the necessary condition of second order \emph{(if a function
$f:E\rightarrow\bkR$ has a local minimum at $x$ and is twice
normally differentiable at this point, then not only $f'(x)=0$ but
also $f''(x)\geq 0$)}. Such conditions are proved by the standard
methods on \cite{Bourbaki}.

Consider the natural extensions of the partial derivatives defined
on (\ref{eq:partDer}) to higher orders.
\begin{theorem}
Let a function $L(t,x,v)$ be twice continuously normally
differentiable on $[a,b]\times E\times E$. If functional
$\textsl{J}[\cdot]$ (\ref{eq:probCVar}) has an extremum at a point
$x\in\Ci$, then $x$ satisfies the Legendre condition
\begin{equation}\label{eq:Legcond}
 \frac{\partial^2 L}{\partial v \partial v}\circ \mu_{x}(t)\geq 0\:\:\:\:\:\forall t\in [a,b].
\end{equation}
\end{theorem}
\begin{proof}
As shortcut define $X\equiv\Ci$, which is a locally convex space.
Recall that $J:X\rightarrow\bkR$,
$\textsl{J}^{\,\prime}_K[\cdot]:X\rightarrow\Fu{X}$ and
$\textsl{J}^{\,\prime\prime}_K[\cdot]:X\rightarrow\Fu{X\times X}$.
Theorem~\ref{teo-orlov} ensures that $\textsl{J}^{\,\prime}_K[x]$
is a K-differentiable mapping at $x$, so it is (locally) normally
differentiable at $x$.
We will show that $\textsl{J}^{\,\prime\prime}_K[x]$ is also
a K-differentiable mapping at $x$.\\
Consider an arbitrary convex compact set $C\subset X$ where $x$ is
a limiting point. The sets $A=\{y\in E : y\in x([a,b]), x\in C\}$
and $B=\{z\in E : z\in \dot{x}([a,b]), \dot{x}\in C\}$ are convex
compacts in $E$. By the definition of normal differentiability, we
have
\begin{gather*}
\textsl{J}^{\,\prime}_K[x+h_2]h_1-\textsl{J}^{\,\prime}_K[x]h_1
=\int_a^b\left[
    \frac{\partial L}{\partial x}(t,x+h_2,\dot{x}+h_2)h_1
     \right.\\
  \left.+\frac{\partial L}{\partial v}(t,x+h_2,\dot{x}+h_2)\dot{h}_1\right]
-\left[\frac{\partial L}{\partial x}(t,x,\dot{x})h_1
     +\frac{\partial L}{\partial v}(t,x,\dot{x})\dot{h}_1 \right]\,dt\\
= \int_a^b \left[ \left(\frac{\partial^2 L}{\partial v\partial v}\circ \mu_{x}\right) \dot{h}_1\dot{h}_2
+\left(\frac{\partial^2 L}{\partial v\partial x}\circ \mu_{x}\right)\dot{h}_2 h_1\right.\\
+\left.\left(\frac{\partial^2 L}{\partial x\partial v}\circ \mu_{x}\right)\dot{h}_1 h_2
+\left(\frac{\partial^2 L}{\partial x\partial x}
\circ \mu_{x}\right)h_1 h_2\right]\, dt
+\int_a^b r_t(h_2,\dot{h}_2)\,dt,
\end{gather*}
where $r_t$ is the residual term of the increments at the point
$(t,x(t),\dot{x}(t))$. Since $L$ is twice continuously normally
differentiable on $C\equiv[a,b]\times A\times B$, we can apply
lemma~\ref{fund:lemma} to $r_t$
\[ \frac{|r_t(h_2(t),\dot{h}_2(t))|}{\|h_2(t)\|_m+\|\dot{h}_2(t)\|_m}\rightarrow 0\mbox{ as }
(h_2(t),\dot{h}_2(t))\rightarrow 0. \] Let $R_x(h_2)=\int_a^b
r_t(h_2(t),\dot{h}_2(t))\,dt$. The uniform convergence and the
above condition imply $\frac{R_x(h_2)}{\|h_2\|_1^m}\rightarrow 0$
as $h_2\rightarrow 0$, which implies the K-differentiability of
$\textsl{J}^{\,\prime\prime}_K[x]$ at $x$.

Hence, the second variation for $\textsl{J}\left[x\right]$
reads as follows, for all $x,h \in\Ci$,
\begin{eqnarray*}\label{eq:secvar}
\textsl{J}''_K\left[x\right](h,h) &=&
\int_a^b \left[ \left(\frac{\partial^2 L}{\partial v\partial v}\circ \mu_{x}(t)\right)
\dot{h}(t)^2+2\left(\frac{\partial^2 L}{\partial v\partial x}
\circ \mu_{x}(t)\right)\dot{h}(t)h(t)\right]\,dt\\
&\ & +\int_a^b \left(\frac{\partial^2 L}{\partial x\partial x}
\circ \mu_{x}(t)\right)h(t)^2\, dt,
\end{eqnarray*}
with the corresponding Jacobi eigenvalue equation
\begin{equation*}
 -\frac{d}{dt}\left(R(t)\dot{h}(t)\right)+P(t)h(t)=\lambda h(t),\:\:\:\:h\in\Ci,
\end{equation*}
where we have defined
\begin{equation*}
 R(t)=\frac{\partial^2 L}{\partial v\partial v}\circ \mu_{x}(t)
 \:\:\:\mbox{ and }\:\:\:
 P(t)=\frac{\partial^2 L}{\partial x\partial x}
\circ \mu_{x}(t)-\frac{d}{dt}\frac{\partial^2 L}{\partial v\partial x}\circ \mu_{x}(t).
\end{equation*}
Now, if we observe that $\frac{d}{dt}h(t)^2=2\dot{h}(t)h(t)$,
then integration by parts of the second variation yields
\begin{equation}\label{eq:secvar2}
  \textsl{J}''_K\left[x\right]h = \int_a^b R(t)\dot{h}(t)^2+P(t)h(t)^2\, dt
  \:\:\:\:\:\:\:\forall x,h \in\Ci.
\end{equation}
The necessary condition of second order implies that
$\textsl{J}''_K\left[x\right](h,h)\geq 0$ for all $h\in\Ci$.
Therefore, (\ref{eq:Legcond}) follows from (\ref{eq:secvar2}).
Namely, if $R(t_0)<0$ for a $t_0\in[a,b]$, then one can always
choose an $h\in\Ci$ having very large $\dot{h}(t_0)$ and small
$h(t_0)$, so that $\textsl{J}''_K\left[x\right](h,h)<0$ holds;
however, this is not possible.
\end{proof}


\section{Invariance and conservation laws}

To obtain a Noether theorem on locally convex spaces,
we will need further regularity
of the solution curves $x$ on the LCS $E$. Similarly to $\Ci$,
we denote by $\Cii$ the space of twice continuous differentiable
mappings $x:[a,b]\rightarrow E$ with a determining system
of semi-norms $\{\|x\|_2^p\}_{p\in P}$ where
\begin{equation*}
\|x\|_2^p=\sup_{a\leq t\leq b}\|x(t)\|_p+\sup_{a\leq t\leq
b}\|\dot{x}(t)\|_p+\sup_{a\leq t\leq b}\|\ddot{x}(t)\|_p.
\end{equation*}
Therefore, the problem of calculus of variations is
\begin{equation}
\label{eq:probCVar2}
\textsl{J}\left[x\right] =
\int_a^b L\left(t,x(t),\dot{x}(t) \right)dt \longrightarrow
\min\, ,
\end{equation}
where $L:[a,b]\times E\times E\rightarrow \bkR$
is twice continuously normally differentiable and $x \in \Cii$.

Let us introduce a local Lie group $h^s$ with generators $T$ e $X$.
A transformation $h$ in the space $\bkR\times E$
is a twice continuously normally differentiable mapping
$h:\bkR\times E  \to \bkR\times E$ with
$h(t,x)=(\bar{t},\bar{x})$ defined by the equations
\begin{equation*}
\bar{t} = h_t\left(t,x\right) \, , \quad
\bar{x} = h_{x}\left(t,x\right)\, ,
\end{equation*}
for $h_t$ and $h_x$ given functions.
The symmetry transformations, which define the invariance of problem
\eqref{eq:probCVar2}, are transformations which depend
on a real parameter $s\in\bkR$. Let $s$ vary continuously in an open interval
$|s| < \varepsilon$, for small $\varepsilon$,
and $h^s$ be a family of transformations defined by
\begin{equation}
\label{eq:ftransformacoes}
\bar{t} = h_t^s(t,x) = h_t\left(t,x,s\right) \, , \quad
\bar{x} = h_x^s(t,x) = h_x\left(t,x,s\right) \, ,
\end{equation}
where $h_t$ and $h_x$ are analytical functions in
$\left[a,b\right] \times E  \times \left(-
\varepsilon,\varepsilon\right)$. The one-parameter family of
transformations $h^s$ is a local Lie group if and only if it
satisfies the \emph{local closure property}; contains the
\emph{identity} (without loss of generality, we assume that the
identity transformation is obtained for parameter $s = 0$); and
\emph{inverse} exist for each $s$ sufficiently small. Since
normally differentiable mappings admit a Taylor formula in the
asymptotic form \cite{Orlov02}, if $h^s$ defined by
\eqref{eq:ftransformacoes} is a local Lie group, then one can
expand $h_t=(t,x,s)$ and $h_{x}=(t,x,s)$ in Taylor series about
$s=0$:
\begin{equation}
\label{eq:tsas0}
\begin{gathered}
\begin{split}
\bar{t} &= h_t\left(t,x,0\right) + \frac{\partial h_t}{\partial s}\left(t,x,0\right) s
+ \frac{1}{2}\frac{\partial ^2 h_t}{\partial s^2}\left(t,x,0\right) s^2  + \cdots \\
&= t + T\left(t,x\right) s + o\left(s\right) \, ,
\end{split}\\
\begin{split}
\bar{x} &= h_{x}\left(t,x,0\right) + \frac{\partial h_{x}}{\partial s}\left(t,x,0\right) s
+ \frac{1}{2}\frac{\partial^2 h_{x}}{\partial s^2}\left(t,x,0\right) s^2  + \cdots \\
&= x + X\left(t,x\right) s + o\left(s\right) \, ,
\end{split}
\end{gathered}
\end{equation}
where the quantities $T(t,x) = \frac{\partial h_t}{\partial s}\left(t,x,0\right)$,
and $X(t,x) = \frac{\partial h_x}{\partial s}\left(t,x,0\right)$
are called the \emph{generators} of $h^s$. A local Lie group $h^s$
induces, in a natural way, a local Lie group $\tilde{h}^s$
in the space $\bkR\{t\}\times E\{x\}\times E\{\dot{x}\}$:
\begin{equation*}
\tilde{h}^s :
\begin{cases}
\bar{t} &= h_t\left(t,x,s\right) \, , \\
\bar{x} &= h_x\left(t,x,s\right) \, , \\
\dot{\bar{x}} &= \frac{d \bar{x}}{d t} = \frac{\frac{\partial h_x}{\partial t}
+ \frac{\partial h_x}{\partial x} \dot{x}}{\frac{\partial
h_t}{\partial t} + \frac{\partial h_t}{\partial x} \dot{x}} \, .
\end{cases}
\end{equation*}
This group is called in the $\bkR^n$-setting the \emph{extended group}. Noticing that
\begin{equation*}
\begin{gathered}
\frac{\partial h_t}{\partial t} = 1 + s \frac{\partial T}{\partial t} + o(s) \, , \quad
\frac{\partial h_t}{\partial x} = s \frac{\partial T}{\partial x} + o(s) \, , \\
\frac{\partial h_x}{\partial t} = s \frac{\partial X}{\partial t} + o(s) \, , \quad
\frac{\partial h_x}{\partial x} = 1 + s \frac{\partial X}{\partial x} + o(s) \, ,
\end{gathered}
\end{equation*}
then
\begin{equation}
\label{eq:expansaoderiv}
\begin{split}
\dot{\bar{x}} &=\frac{s\frac{\partial X}{\partial t}
+ \left( 1 + s\frac{\partial X}{\partial x}\right) \dot{x} + o(s)}{1
+ s \frac{\partial T}{\partial t} + s\frac{\partial T}{\partial x}\dot {x}
+ o(s)} = \frac{\dot{x} + s X' + o(s)}{1 + s T' + o(s)} \\
&= \dot{x} + \left( X' - \dot{x} T'\right)s + o(s)\\
&= \dot{x} + V s + o(s)\, ,
\end{split}
\end{equation}
with $V(t,x,\dot{x}) = X'(t,x) - \dot{x}T'(t,x)$
the generator associated with the derivative.

The integral functional $\textsl{J}[\cdot]$ of the fundamental problem
of the calculus of variations \eqref{eq:probCVar2},
\begin{equation*}
\textsl{J}\left[x\right] =
\int_a^b L\left(t,x(t),\dot{x}(t)\right)dt \, ,
\end{equation*}
is said to be \emph{invariant} under a local Lie group $h^s$ if, and only if,
\begin{equation}
\label{eq:definv1}
\begin{split}
\int_{t_1}^{t_2} L\left(t,x(t),\frac{dx}{dt}(t)\right)dt + o(s)
&=\int_{\bar{t}_1}^{\bar{t}_2} L\left(\bar{t},\bar{x}\left(\bar{t}\right),
\frac{d\bar{x}}{d\bar{t}}\left(\bar{t}\right)\right)d\bar{t}\\
&= \int_{t_1 }^{t_2 } L\left( h_t^s ,h_{x}^s ,\frac{\frac{\partial h_x^s}{\partial t}
+ \dot{x}\frac{\partial h_x^s}{\partial x}}{\frac{\partial h_t^s}{\partial t}
+ \dot{x}\frac{\partial h_t^s}{\partial x}}\right)\frac{d h_t^s}{dt}dt \, ,
\end{split}
\end{equation}
where $\bar{t}_1 = h_t(t_1,x(t_1),s)$, $\bar{t}_2 = h_t(t_2,x(t_2),s)$,
$h_t^s=h_{t}(t,x(t),s)$, $h_x^s=h_x(t,x(t),s)$,
and \eqref{eq:definv1} is verified
for all $|s| < \varepsilon$, for all $x \in C^2([a,b];E)$,
and for all $[t_{1},t_{2}] \subseteq [a,b]$.
Condition \eqref{eq:definv1} is equivalent to
\begin{multline}
\label{eq:definv1b}
\left.\frac{d}{ds}\int_{\bar{t}_1}^{\bar{t}_2} L\left(\bar{t},\bar{x}\left(\bar{t}\right),
\frac{d\bar{x}}{d\bar{t}}\left(\bar{t}\right)\right)d\bar{t}\right|_{s=0} = 0 \\
\Leftrightarrow \left.\frac{d}{ds}
\int_{t_1}^{t_2} L\left( h_t^s ,h_x^s,\frac{\frac{\partial h_x^s}{\partial t}
+ \dot{x}\frac{\partial h_x^s}{\partial x}}{\frac{\partial h_t^s}{\partial t}
+ \dot{x}\frac{\partial h_t^s}{\partial x}}\right)\frac{d h_t^s}{dt}dt\right|_{s=0} = 0 \, .
\end{multline}
The requirement that \eqref{eq:definv1b} hold for every subinterval $[t_1,t_2]$
of $[a,b]$ permits to remove the integral from consideration,
and to put the focus on the Lagrangian $L(\cdot,\cdot,\cdot)$.

\begin{definition}[Invariance/symmetry]
The fundamental problem of the calculus of variations on locally convex spaces
\eqref{eq:probCVar2} is said to be \emph{invariant} under the local Lie group
$h^s$ if, and only if,
\begin{equation}
\label{eq:eqsintinv}
\left.\frac{d}{ds} \left\{
 L\left( h_t^s ,h_x^s,\frac{\frac{\partial h_x^s}{\partial t}
+ \dot{x}\frac{\partial h_x^s}{\partial x}}{\frac{\partial h_t^s}{\partial t}
+ \dot{x}\frac{\partial h_t^s}{\partial x}}\right)\frac{d h_t^s}{dt}\right\}\right|_{s=0} = 0 \, .
\end{equation}
We then say that $h^s$ is a \emph{symmetry} for the problem.
\end{definition}

\begin{theorem}[Necessary and sufficient condition of invariance]
The fundamental problem of the calculus of variations \eqref{eq:probCVar2}
is invariant under a local Lie group
\eqref{eq:ftransformacoes} with generators $T$ and $X$ if, and only if,
\begin{multline}
\label{eq:condnsinv}
(L^t_{(x(t),\dot{x}(t))})'(t) T(t,x(t)) + (L^x_{(t,\dot{x}(t))})'(x(t)) X(t,x(t))\\
+ (L^v_{(t,x(t))})'(\dot{x}(t))\bigl(X'(t,x(t)) - \dot{x}(t)T'(t,x(t))\bigr)
+ L \circ \mu_{x}(t) T'(t,x(t)) = 0 \, .
\end{multline}
\end{theorem}
\begin{proof}
Carrying out the differentiation of equation \eqref{eq:eqsintinv} we obtain:
\begin{equation}
\label{eq:quaseProv}
\left. L\frac{d}{ds}\left(\frac{d\bar{t}}{dt}\right)\right|_{s = 0}
+ \left.\frac{d}{ds} L\left(\bar{t},\bar{x},\frac{d\bar{x}}{d\bar{t}}\right)\right|_{s = 0}  = 0 \, .
\end{equation}
Recalling \eqref{eq:partDer}, and
having in mind that by \eqref{eq:tsas0} and \eqref{eq:expansaoderiv}
\begin{equation*}
\begin{split}
\left.\frac{d}{ds}\left(\frac{d\bar{t}}{dt}\right)\right|_{s = 0}
&= \left.\frac{d}{ds}\left(\frac{d}{dt}\left(t + sT + o(s)\right)\right)\right|_{s = 0} = T' \, , \\
\left.\frac{d}{ds} L\left(\bar{t},\bar{x},\frac{d\bar{x}}{d\bar{t}}\right)\right|_{s = 0}
&= \frac{\partial L}{\partial t} T + \frac{\partial L}{\partial x} X
+ \frac{\partial L}{\partial v}\left(X' - \dot{x}T'\right) \, ,
\end{split}
\end{equation*}
we obtain from \eqref{eq:quaseProv} the intended conclusion.
\end{proof}

\begin{theorem}[Noether's Theorem on Locally Convex Spaces]
\label{thm:MR:NT}
If
\begin{equation*}
\begin{gathered}
\textsl{J}\left[x\right] = \int_a^b L\left(t,x(t),\dot{x}(t)\right)dt
\end{gathered}
\end{equation*}
is invariant under a local Lie group $h^s$ with
generators $T$ and $X$, then
\begin{multline}
\label{eq:prob2'}
\left[ L \circ \mu_{x}(t) - \dot{x}(t)
(L^v_{(t,x(t))})'(\dot{x}(t)) \right] T(t,x(t))\\
+ (L^v_{(t,x(t))})'(\dot{x}(t)) X(t,x(t)) = \text{constant}
\end{multline}
$\forall t \in [a,b]$, and along all the solutions $x$
of the Euler-Lagrange equation \eqref{eq:EL:LCS}.
\end{theorem}

\begin{remark}
If we introduce the Hamiltonian function $H(\cdot,\cdot,\cdot)$ by
\begin{equation*}
H\left(t,x,\dot{x}\right) = - L\left(t,x,\dot{x}\right)
+ \dot{x} \frac{\partial L}{\partial v}\left(t,x,\dot{x}\right) \, ,
\end{equation*}
with $\frac{\partial L}{\partial v}$ as in \eqref{eq:partDer},
then the conservation law \eqref{eq:prob2'} can be written in the form
\begin{equation*}
\left[\frac{\partial L}{\partial v} \circ \mu_{x}(t)\right] X\left(t,x(t)\right)
- \left[H \circ \mu_{x}(t)\right] T\left(t,x(t)\right) = \text{constant} \, .
\end{equation*}
\end{remark}

\begin{proof}
Direct calculations show that:
\begin{gather}
\frac{d}{dt}\left( \frac{\partial L}{\partial v} X \right)
- X\frac{d}{{dt}}\frac{\partial L}{\partial v} =
\frac{d} {{dt}}\frac{\partial L}{\partial v} X
+ \frac{\partial L}{\partial v} X'
- X\frac{d}{{dt}}\frac{\partial L}{\partial v}
= \frac{\partial L}{\partial v} X' \, , \label{eq:prob1} \\
\begin{split}
\frac{{dL}}{{dt}} T - \frac{\partial L}{\partial x} \dot{x} T
- \frac{\partial L}{\partial v} \ddot{x} T &=\left(
{\frac{{\partial L}} {{\partial t}} + \frac{{\partial L}}
{{\partial x}}\dot x + \frac{{\partial L}} {{\partial v}}\ddot x} \right) T
- \frac{\partial L}{\partial x} \dot{x} T
- \frac{\partial L}{\partial v} \ddot{x} T \\
&= \frac{{\partial L}}{{\partial t}} T \, ,
\end{split} \label{eq:prob2} \\
\begin{split}
\frac{d}{{dt}}\left( {\frac{\partial L}{\partial v} \dot{x} T} \right)
- \dot{x} T\frac{d}{{dt}}\frac{\partial L}{\partial v}
&=  \frac{d} {{dt}}\frac{\partial L}{\partial v}
\dot{x} T + \frac{\partial L}{\partial v}
\frac{d} {{dt}}\left( \dot{x} T \right)
- \dot x T\frac{d}{{dt}}\frac{\partial L}{\partial v}\\
&= \frac{\partial L}{\partial v} \left( {\ddot x T + \dot x T'} \right) \, .
\end{split} \label{eq:prob3}
\end{gather}
Substituting \eqref{eq:prob2} in the necessary and sufficient
condition of invariance \eqref{eq:condnsinv} we obtain that
\begin{equation}
\label{eq:ax1}
\frac{d L}{{dt}} T - \frac{{\partial L}}
{{\partial x}}\dot x T - \frac{{\partial L}}
{{\partial v}}\ddot x T + \frac{{\partial L}}
{{\partial x}} X
+ \frac{{\partial L}}{{\partial v}} X' - \frac{{\partial L}}
{{\partial v}}\left( {\dot x T'} \right) + L T' = 0 \, ;
\end{equation}
substituting \eqref{eq:prob1} in \eqref{eq:ax1} we get
\begin{multline}
\label{eq:ax2}
\frac{d L}{{dt}} T - \frac{{\partial L}}
{{\partial x}}\dot x T - \frac{{\partial L}}
{{\partial v}}\ddot x T +
\frac{{\partial L}}
{{\partial x}} X +
\frac{d}{dt}\left(\frac{{\partial L}}{{\partial v}} X\right)
- X\frac{d}{{dt}}\frac{{\partial L}}
{{\partial v}} - \frac{{\partial L}}
{{\partial v}}\left( {\dot x T'} \right) + L T' = 0 \, .
\end{multline}
Finally, using \eqref{eq:prob3} in the last equation
\eqref{eq:ax2} one obtains
\begin{multline*}
\frac{dL}{{dt}} T - \frac{{\partial L}}
{{\partial x}}\dot x T - \frac{{\partial L}}
{{\partial v}}\ddot x T + \frac{{\partial L}}
{{\partial x}} X + \frac{d}
{{dt}}\left( {\frac{{\partial L}}
{{\partial v}} X} \right) \\
- X\frac{d}
{{dt}}\frac{{\partial L}}
{{\partial v}} + \frac{{\partial L}}
{{\partial v}}\ddot x T - \frac{d}
{{dt}}\left( {\frac{{\partial L}}
{{\partial v}}\dot x T} \right) + \dot x T\frac{d}
{{dt}}\frac{{\partial L}}
{{\partial v}} +  L T' = 0  \, ,
\end{multline*}
which, by simplification, takes form
\begin{multline*}
\frac{dL}{{dt}} T
+ \frac{{\partial L}}{{\partial x}}\left( {X - \dot x T} \right)
+ \frac{d}{{dt}}\left( {\frac{{\partial L}}
{{\partial v}} X} \right) -
\frac{d} {{dt}} \left(\frac{{\partial L}}
{{\partial v}}{\dot x T} \right) - \frac{d}
{{dt}}\frac{{\partial L}}
{{\partial v}}\left( {X - \dot x T} \right) + L T' = 0  \\
\Leftrightarrow \frac{d}
{{dt}}\left( {L T + \frac{{\partial L}}
{{\partial v}} X - \frac{{\partial L}}
{{\partial v}}\dot x T} \right)
+ \left( {X - \dot x T} \right)\left( {\frac{{\partial L}}
{{\partial x}} - \frac{d}
{{dt}}\frac{{\partial L}}
{{\partial v}}} \right) = 0  \, .
\end{multline*}
By definition, along the solutions of the Euler-Lagrange equations
\eqref{eq:EL:LCS} $\frac{\partial L}{\partial x}
- \frac{d}{dt}\frac{\partial L}{\partial v} = 0$,
and one gets the desired conclusion:
\begin{equation*}
\frac{d}{dt}\left(L T + \frac{\partial L}{\partial v} X
- \frac{\partial L}{\partial v}\dot{x} T \right) = 0  \, .
\end{equation*}
\end{proof}

The theory of the calculus of variations on locally convex spaces
is under development. Much remains to be done, in particular, in
the vectorial setting. Even for the scalar case, it would be
interesting, for example, to have a version of the DuBois-Reymond
necessary condition on locally convex spaces. With such condition
one can try to prove Theorem~\ref{thm:MR:NT} for more general
classes of admissible functions, following the scheme in
\cite{Torres04}.


\section*{Acknowledgment}

The authors acknowledge the financial support of the FCT/FEDER
project POCTI/MAT/41683/2001 and the Research Unit CEOC. The first
author also acknowledge the financial support of the FCT/FEDER
project POCI\-/MAT\-/55524\-/2004.



\end{document}